\documentclass{article}

\input{amssym}
\begin{document}
\title{{\bf Symmetry of Quadratic Homogeneous Differential Systems}}
\date{}
\author{Mehdi Nadjafikhah\thanks{Department of
Mathematics, Iran University of Science and Technology, Narmak-16,
Tehran, Iran. e-mail: m\_nadjafikhah@iust.ac.ir} \and Ali
Mahdipour-Shirayeh\thanks{mahdipour@iust.ac.ir}}
\maketitle
\begin{abstract} In this paper, the symmetry group of a differential system of
$n$ quadratic homogeneous first order ODEs of $n$ variables is
studied. For this purpose, we consider the action of both point
and contact transformations to signify the corresponding Lie
algebras. We also find the independent differential invariants of
these actions.
\end{abstract}
\noindent {\bf A.M.S. 2000 Subject Classification}: 76M60, 35E20, 53A55.
\section{Introduction}

A quadratic homogeneous differential system (QHDS) is a
generalization of second degree polynomial system of differential
equations that is rooted in natural science, and also in
mathematics by virtue of Hilbert's sixteen problem of his famous
list of problems (See, e.g., \cite{CJ} and references therein).
They are studied in some different view point of their
applications. For example, in physical sciences, in van der Pol
oscillator as an important example of qualitative theory of
ordinary differential equations (ODEs), in mathematical echology
and in particular in Volterra-Lotka equations, and many other
applications in astrophysics and fluid mechanics \cite{CJ}.\par
There are many aspects for studying a QHDS.  For instance, in
\cite{Ne} geometric classification of trajectories of a given two
dimensional QHDS has studied for determining the invariant lines
through the origin of the QHDS and also its location in the plane
of parameters for classification of the geometry of its
trajectories. In \cite{Ts}, the necessary conditions for the
existence of polynomial first integrals and a polynomial symmetry
field of a n-dimensional QHDS have given. \par At the present
work, the symmetry group of a QHDS of every finite dimension is
studied. We specify symmetry groups of both point and contact
transformation groups.\\

Assume that we have a differential system ${\cal E}$ of $n$
quadratic homogeneous first order ODEs of $n$ variables, in which
each of them is a function of $t$:
\begin{eqnarray}
{\cal E}:\;\;\dot{x}_i = F_i(x_1,\cdots,x_n), \hspace{1cm}F_i
=\sum_{j,k}\,a^i_{jk}\,x_j\,x_k,\label{eq:1}
\end{eqnarray}
where the coefficients, $a^i_{j,k}$ for $1\leq i, j, k\leq n$, are
arbitrary real constants. The solution space of differential
system (\ref{eq:1}) is an immersed submanifold of the first order
jet space, ${\rm J}^1({\Bbb R},{\Bbb R}^n)$ of dimension $2n+1$
and local coordinates $(t,{\bf x},{\bf p}):=
(t,x_1,\cdots,x_n,p_1,\cdots,p_n)$, when each $p_j$ is just
$\dot{x}_j$.\par
According to \cite{Ol}, the action of contact transformations
group is introduced on the jet space ${\rm J}^1({\Bbb R},{\Bbb
R}^n)$, and any element of this group is as the following
transformation:
\begin{eqnarray}
T=\chi(t,{\bf x},{\bf p}),\hspace{0.4cm}X_i=\psi_i(t,{\bf x},{\bf
p}),\hspace{0.4cm}P_j=\pi_j(t,{\bf x},{\bf p}),\label{eq:2}
\end{eqnarray}
where $1\leq i, j\leq n$. An important special case of contact
transformation group is point transformation group, where the
action reduces to the following change of coordinates
\begin{eqnarray*}
T=\chi(t,{\bf x}),\hspace{0.4cm}X_i=\psi_i({\bf x},{\bf p}),
\end{eqnarray*}
when the variables come from the jet space $J^0({\Bbb R},{\Bbb
R}^n)$.\par
\paragraph{Remark 1.1} {\it In general, the point and contact symmetry groups of a
geometric object are not necessarily equal. A jet space $J^k({\Bbb
R}^m,{\Bbb R}^n)$ of order $k$ can be imbedded in other jet spaces
$J^l({\Bbb R}^m,{\Bbb R}^n)$ of order $l>k$ as a submanifold. Also
a contact transformation of $J^l({\Bbb R}^m,{\Bbb R}^n)$ can be
restricted to $J^k({\Bbb R}^m,{\Bbb R}^n)$ and then acts as a
point transformation. Thus, the point transformation group of an
geometric object is a Lie group of the contact transformation
group and in the special case of symmetries, we conclude that the
point symmetry group is a Lie subgroup of the contact symmetry group.}\\

In the next section, we will see that point symmetry group of the
differential system ${\cal E}$ is of dimension 1, where its
contact symmetry group has infinite dimension.
\section{Point Symmetry of a QHDS}
As is indicated above, an important question of a quadratic
homogeneous system is about its symmetry group. There are various
methods for signifying symmetry group based on Lie's method of
infinitesimal generators, that some of them are exist in
\cite{Ha}, \cite{Ib}, \cite{KLR} and \cite{Ol}. In this work, we
will follow the method of \cite{Ol} for finding point and contact
symmetries.\\

In jet space $J^1({\Bbb R},{\Bbb R}^n)$ the system of equations
(\ref{eq:1}) will be as follows
\begin{eqnarray}
p_i=\sum_{j,k}\,a^i_{jk}\,x_j\,x_k,\hspace{0.8cm} 1\leq i, j,
k\leq n.\label{eq:3}
\end{eqnarray}
For finding symmetry group of system ${\cal E}$, we consider the
general form of the infinitesimal generators provided from a point
transformation as follows:
\begin{eqnarray}
v=T\frac{\partial}{\partial t}+
\sum_{i=1}^n\,K_i\frac{\partial}{\partial x_i},\label{eq:4}
\end{eqnarray}
in which, $T$ and $K_i$s are smooth functions of variables $t,
x_1, \cdots, x_n$.\par
The fist order prolongation of $v$ will be
then as the following expression:
\begin{eqnarray} v^{(1)}= v +
\sum_{i=1}^n\,K_i^t\frac{\partial}{\partial p_i},\label{eq:5}
\end{eqnarray}
where for each $1\leq i\leq n$ we have
\begin{eqnarray}
K_i^t=K_{i,t}+\sum_j K_{i,x_j}-p_i(T_t+\sum_j
T_{i,x_j}).\label{eq:6}
\end{eqnarray}
According to \cite{Ol}, $v^{(1)}$ as a prolongation of $v$ is a
symmetry of system ${\cal E}$, if and only if it satisfies in
relation $v^{(1)}[{\cal E}]=0$. By effecting $v^{(1)}$ on the
equations of the system we conclude the following system of
equations for $i=1,\cdots, n$:
\begin{eqnarray}
&&K_1\Big(\sum_k\,(a^i_{1k}+a^i_{k1})x_k\Big) + \cdots +
K_n\Big(\sum_k\,(a^i_{nk}+a^i_{kn})x_k\Big)\nonumber\\
&&- K_{i,t} - p_1\,K_{i,x_1} - p_2\,K_{i,x_2} - \cdots - p_n\,K_{i,x_n}\label{eq:7}\\
&&+ p_i(T_t + p_1\, T_{x_1} + \cdots + p_n\, T_{x_n})=0. \nonumber
\end{eqnarray}
In the last equation, since $t, x_i, p_ j$ ($1\leq i,j\leq n$) are
arbitrary and functions $K_i$ and $N_j$ only depend on $t,x_i$, so
equations (\ref{eq:7}) will be satisfied if and only if we have
the following system of equations
\begin{eqnarray}
&&\hspace{-2cm} K_1\Big(\sum_k\,(a^i_{1k}+a^i_{k1})x_k\Big) +
\cdots + K_n\Big(\sum_k\,(a^i_{nk}+a^i_{kn})x_k\Big)- K_{i,t}=0, \label{eq:8}\\
&&\hspace{-2cm} K_{i,x_1}=0,\hspace{0.1cm}
K_{i,x_2}=0,\hspace{0.1cm} \cdots, \hspace{0.1cm}K_{i,x_{i-1}}=0,
\hspace{0.1cm} K_{i,x_{i+1}}=0,\hspace{0.1cm}\cdots,\hspace{0.1cm} K_{i,x_n}=0, \label{eq:9}\\
&&\hspace{-2cm} T_t - K_{i,x_i}=0,\hspace{0.4cm} T_{x_1}=0,
\hspace{0.4cm} \cdots,\hspace{0.4cm}T_{x_n}=0. \label{eq:10}
\end{eqnarray}
for every $1\leq i\leq n$. At first, from Eqs.
(\ref{eq:9})--(\ref{eq:10}) we understand that $T$ depends only to
$t$ and for each $i$, $K_i$ is a function of $t$ and $x_i$. For a
fixed $1\leq i\leq n$, by replacing the Eq. $T_t = K_{i,x_i}$ in
(\ref{eq:8}) and then differentiating with respect to $x_i$, we
will find the following expression
\begin{eqnarray*}
T_t\Big(\sum_k\,(a^i_{ik}+a^i_{ki})x_k\Big)=T_{tt},
\end{eqnarray*}
By solving the last equation in respect to $t$, we have the
following solution
\begin{eqnarray*}
T(t)=c_1\,\exp\Big(t\,\sum_k\,(a^i_{ik}+a^i_{ki})x_k\Big)+c_2,
\end{eqnarray*}
for constants $c_1$ and $c_2$. Since $T$ is just a function of
$t$, so $c_1$ must be zero, and therefore
$$T(t,x_1,\cdots,x_n)=c,$$
for arbitrary constant $c$. Then for each $i$, we see that
$K_{i,x_i}=0$ and so is dependent to $t$ only. But by solving Eqs.
(\ref{eq:8}) with respect to $K_i(t)$ we obtain the below
solutions when $i$ varies between 1 and $n$:
\begin{eqnarray}
K_i(t)&=&\exp\Big(t\,\sum_k\,(a^i_{ik}+a^i_{ki})x_k\Big)\nonumber \\
&&\left(\int\,\exp\Big(-t\,\sum_k\,(a^i_{ik}+a^i_{ki})x_k\Big)\,F\,dt+c_3\right),\label{eq:11}
\end{eqnarray}
where we assumed that $F=\sum_{j\neq
i}\,K_j(t)\Big(\sum_k\,(a^i_{jk}+a^i_{kj})x_k\Big)$ and $c_3$ be
arbitrary constant. Eqs. (\ref{eq:11}) and the fact that all of
$K_i$s are just depend to $t$ implies that
\begin{eqnarray*}
K_i(t,x_1,\cdots,x_n)=0\hspace{1cm} \mbox{for $1\leq i\leq n$}.
\end{eqnarray*}
when $c$ is an arbitrary constant. The symmetry group of the
general infinitesimal generators, that provided as $v=c
\frac{\partial}{\partial t}$ will be just the translation of $t$
by a constant coefficient of the parameter $s$ of the 1-parameter
subgroup:
\begin{eqnarray*}
(t,x_1,\cdots,x_n)\mapsto(c\,s + t,x_1,\cdots,x_n)
\end{eqnarray*}
\paragraph{Theorem 2.1} {\it Set of all infinitesimal generators of the point symmetry group
of the system (\ref{eq:1}) is a 1-dimensional Lie algebra with the
base $\left\{\frac{\partial}{\partial t}\right\}$. Therefore, the
symmetry group of the system (\ref{eq:1}) is a 1-dimensional Lie
group of time translations.}\\

According to theorem 2.74 of \cite{Ol}, the invariants
$u=I(t,x_1,x_2,\cdots,x_n)$ of one--parameter group with
infinitesimal generators of the form (\ref{eq:4}) satisfy the
linear, homogeneous partial differential equations of first order:
\begin{eqnarray*}
v[I]=0.
\end{eqnarray*}
The solutions of the last equation, are found by the method of
characteristics (See \cite{Ol} and \cite{Ib} for details). So we
can replace the last equation by the following characteristic
system of ordinary differential equations
\begin{eqnarray}
\frac{dt}{T}=\frac{dx_1}{K_1}=\frac{dx_2}{K_2}=\cdots=\frac{dx_n}{K_n}.\label{eq:12}
\end{eqnarray}
By solving the equations (\ref{eq:12}) of the differential
generator $v=c \frac{\partial}{\partial t}$ , we (locally) find
the following general solutions
\begin{eqnarray}
\hspace{-0.5cm}I_1(t,x_1,x_2,\cdots,x_n) = x_1=d_1,
\hspace{0.4cm}\cdots,\hspace{0.4cm} I_n(t,x_1,x_2,\cdots,x_n) =
x_n=d_n, \label{eq:13}
\end{eqnarray}
for constants $d_i$ when $i=1,\cdots,n$.
\paragraph{Theorem 2.2} {\it The independent first integrals
of the characteristic system of the infinitesimal generator $v=c \frac{\partial}{\partial t}$
are the derived invariants (\ref{eq:13}).}\\
\section{Contact Symmetry of a QHDS}
On the other hand, we consider that the infinitesimal generators
comes from contact transformations, that is, the contact
transformation group is acting on the jet space $J^1$. Hence, we
may suppose general form of a infinitesimal generator as
\begin{eqnarray}
v=T\frac{\partial}{\partial t}+
\sum_{i=1}^n\,K_i\frac{\partial}{\partial x_i} +
\sum_{j=1}^n\,P_j\frac{\partial}{\partial p_j}, \label{eq:14}
\end{eqnarray}
for arbitrary smooth functions $T,K_i,P_j$ of variables
$t,x_i,p_j$.\par Since this expression depends to variables of
1-jet space and so is a infinitesimal generator of the jet space,
so it dose not need to be prolonged. Therefore, according to
\cite{Ol} $v$ is a symmetry of ${\cal E}$, if and only if, it
satisfies in $v[{\cal E}]=0$. After acting $v$ on the system we
obtain the relation
\begin{eqnarray*}
K_1\sum_k\,(a^i_{1k}+a^i_{k1})x_k + \cdots +
K_n\sum_k\,(a^i_{nk}+a^i_{kn})x_k- P_i=0
\end{eqnarray*}
where $i$ varies between 1 and $n$. Thus the general form of each
symmetry of the system will be as follows
\begin{eqnarray}
v &=& T\frac{\partial}{\partial t}+
\sum_{i=1}^n\,\Big\{K_i\frac{\partial}{\partial x_i}
+\Big(\sum_{j,k}\,K_j(a^i_{jk}+a^i_{kj})x_k \Big)
\frac{\partial}{\partial p_i}\Big\}, \label{eq:15}
\end{eqnarray}
in which $T,K_1,K_2,\cdots,K_n$ are arbitrary functions of
$t,x_1,x_2,\cdots,x_n$. One may divides $v$ to the following
vector fields
\begin{eqnarray}
v_1 &=& T\frac{\partial}{\partial t} \nonumber\\
v_{2} &=& K_1\Big(\frac{\partial}{\partial x_1}+
\sum_{j,k}\,(a^j_{1k}+a^j_{k1})x_k\,\frac{\partial}{\partial
p_j}\Big)  \nonumber\\
&\vdots&   \label{eq:16}\\
v_{n+1} &=& K_n\Big(\frac{\partial}{\partial x_n}+
\sum_{j,k}\,(a^j_{nk}+a^j_{kn})x_k\,\frac{\partial}{\partial
p_j}\Big), \nonumber
\end{eqnarray}
Lie bracket of each two of $v_1,\cdots, v_{n+1}$ is defined to be
the commutator $[v_i,v_j]:=v_i\,v_j-v_j\,v_i$ ($1\leq i,j\leq
n+1$). It is not hard to see that the commutator of any two of
$v_1,\cdots, v_{n+1}$ is a linear combination of them, and hence
$v_1,\cdots, v_{n+1}$ generate a Lie algebra. For instance for
$1\leq\alpha,\beta\leq n$ if we assume that
\begin{eqnarray*}
v_{\alpha+1}&=& K_{\alpha}\Big(\frac{\partial}{\partial
x_{\alpha}}+
\sum_{j,k}\,(a^j_{{\alpha}k}+a^j_{k{\alpha}})x_k\,\frac{\partial}{\partial
p_j}\Big),\\
v_{\beta+1}&=& K_{\beta}\Big(\frac{\partial}{\partial x_{\beta}}+
\sum_{l,m}\,(a^l_{{\beta}m}+a^l_{m{\beta}})x_m\,\frac{\partial}{\partial
p_l}\Big),
\end{eqnarray*}
then we have
\begin{eqnarray*}
&&[v_{\alpha+1},v_{\beta+1}] =\{K_{\alpha}\,K_{\beta,x_{\alpha}}
+K_{\alpha}\sum_{j,k}\,K_{\beta,p_j}(a^j_{{\alpha}k}+a^j_{k{\alpha}})x_k\}\frac{\partial}{\partial
x_{\beta}}\\
&&\hspace{1cm}-\{K_{\beta}\,K_{\alpha,x_{\beta}}+K_{\beta}\sum_{j,k}\,K_{\alpha,p_j}(a^j_{{\beta}k}
+a^j_{k{\beta}})x_k\}\frac{\partial}{\partial x_{\alpha}}\\
&&\hspace{1cm} +
\sum_j\,\Big\{K_{\alpha}\,K_{\beta,x_{\alpha}}\sum_k\,(a^j_{{\beta}k}+a^j_{k{\beta}})x_k\\
&&\hspace{1cm} -K_{\beta}\,K_{\alpha,x_{\beta}}\sum_k\,(a^j_{{\alpha}k}+a^j_{k{\alpha}})x_k\\
&&\hspace{1cm}+\Big(K_{\alpha}\sum_{l,m}\,K_{\beta,p_l}(a^j_{{\alpha}m}+a^j_{m{\alpha}})x_m\Big)\Big(\sum_{k}
\,(a^j_{{\beta}k}+a^j_{k{\beta}})x_k\Big)\\
&&\hspace{1cm}-K_{\beta}\Big(\sum_{l,m}\,K_{\alpha,p_l}(a^j_{m{\beta}}+a^j_{{\beta}m})x_m\Big)\Big(\sum_{k}
\,(a^j_{{\alpha}k}+a^j_{k{\alpha}})x_k\Big)\Big\}\,
\frac{\partial}{\partial p_j}.
\end{eqnarray*}
One can compute the table of Lie symmetry algebra of commutators
arising from infinitesimal operators (\ref{eq:16}) that is given
in table 1.\\
\begin{eqnarray*}
 \begin{array}{|l|ccccc|}
  \hline
        & v_1         & v_2         & \cdots     & v_n         & v_{n+1}    \\ \hline
  v_1   & 0           & v_2-v_1     & \cdots     & v_n-v_1     & v_{n+1}-v_1      \\[1mm]
  v_2   & v_1-v_2     & 0           & \cdots     & v_n-v_2     & v_{n+1}-v_2    \\[1mm]
 \vdots &\vdots       &\vdots       &\ddots      & \vdots      & \vdots    \\[1mm]
  v_n   & v_n-v_1     &v_2-v_n      & \cdots     & 0           &v_{n+1}-v_n    \\[1mm]
 v_{n+1}& v_1-v_{n+1} &v_2-v_{n+1}  &\cdots      &v_n-v_{n+1}  & 0      \\[1mm]
  \hline
 \end{array} \\
 \centering{\mbox{\small Table 1. The commutators table of Lie algebra of equation
 (\ref{eq:14})}
 }\label{tabel:1} \nonumber
\end{eqnarray*}
Therefore we have the following theorem:
\paragraph{Theorem 3.1} {\it The symmetry group of the system ${\cal
E}$ of quadratic homogeneous differential equations, (\ref{eq:1}),
under the action of contact transformations is an infinite
dimensional Lie algebra generated by infinitesimal operators
(\ref{eq:16}).}\\

By referring again to theorem 2.74 of \cite{Ol}, and considering
the infinitesimal operators of the form (\ref{eq:15}), the
invariants $I(t,{\bf x},{\bf p})$ of their one-parameter group,
are signified by the characteristic system of ordinary
differential equations as follows ($i,j,k=1,\cdots,n$)
\begin{eqnarray}
\frac{dt}{T}=\frac{dx_i}{K_i}=\frac{dp_j}{K_1\sum_k\,(a^j_{1k}+a^j_{k1})x_k
+\cdots+K_n\sum_k\,(a^j_{nk}+a^j_{kn})x_k}\label{eq:17}
\end{eqnarray}
in which determines $2n$ invariant functions. By solving these
equations we find that
\begin{eqnarray}
I_1(t,{\bf x},{\bf p}) &=& \int(K_1\,dt-T\,dk_1)=d_1, \nonumber\\
&\vdots& \nonumber \\
I_n(t,{\bf x},{\bf p}) &=& \int(K_n\,dt-T\,dk_n)=d_n, \label{eq:18}\\
I_{n+1}(t,{\bf x},{\bf p}) &=&
\int\Big(\sum_{j,k}\,K_j(a^1_{jk}+a^1_{kj})x_k\Big)dt - \int
T\,dp_1 = d_{n+1}\nonumber\\
&\vdots& \nonumber\\
I_{2n}(t,{\bf x},{\bf p}) &=&
\int\Big(\sum_{j,k}\,K_j(a^n_{jk}+a^n_{kj})x_k\Big)dt - \int
T\,dp_n = d_{2n} \nonumber
\end{eqnarray}
\paragraph{Theorem 3.2} {\it The derived invariants (\ref{eq:18}) as independent first integrals
of the characteristic system of the infinitesimal generator
(\ref{eq:15}), provide the general solution
$$S(t,{\bf x},{\bf p}):=\varphi(I_1(t,{\bf x},{\bf p}),I_2(t,{\bf x},{\bf p}),\cdots,I_{2n}(t,{\bf x},{\bf p})),$$ with an arbitrary function
$\varphi$, which satisfies in the equation $v[\varphi]=0$ and
therefore in
$$T\,\frac{\partial S}{\partial t}+\sum_{i=1}^n \Big(K_i\,\frac{\partial S}{\partial k_i}
+\Big(\sum_{j,k}\,K_j(a^i_{jk}+a^i_{kj})x_k\Big)\,\frac{\partial
S}{\partial p_i}\Big)=0.$$ The solution of the last equation is
defined implicitly by $S(t,{\bf x},{\bf p})=0$. If $\frac{\partial
S}{\partial t}\neq 0$, then the solution can be
written explicitly by $t=\gamma({\bf x},{\bf p})$.}\\
\paragraph{Example 3.3} Clearly, by assuming $T=1$ and $K_i=0$ we
return back to a point symmetry, in which as we saw in previous
section, the symmetry group provides by translating time and
fixing other variables.
\paragraph{Example 3.4} If we suppose $T=1$ and $K_i=1$ for some $i$ and $K_j=0$ for $j\neq
i$, then we have the following expression of infinitesimal
generators
\begin{eqnarray}
v = \frac{\partial}{\partial t}+\frac{\partial}{\partial x_i}+
\Big(\sum_{j,k}\,(a^j_{ik}+a^j_{ki})x_k\Big)\frac{\partial}{\partial
p_j}. \label{eq:19}
\end{eqnarray}
Its one--parameter group then will be as follows
\begin{eqnarray*}
(t,{\bf x},{\bf p})& \mapsto&\Big(t+s,x_1,\cdots, x_i+s, \cdots,
x_n, p_m +
a^m_{ii}s^2+s\,\sum_k\,(a^m_{ik}+a^m_{ki})x_k\Big)\label{eq:20}
\end{eqnarray*}
where $s$ is the parameter of the flow and $m$ changes over
$1,\cdots,n$. If the point $(t,{\bf x},{\bf p})$ be fixed, then
the flow of $v$, by restricting the jet space coordinates to the
$t,x_i$ and $p_j$ axis, for a $1\leq j\leq n$ and the fixed value
$i$, is a space parabolic.
\paragraph{Example 3.5} Let $T=t$, $K_i=1$ for some index $i$ and $K_j=0$ for other indices $j\neq i$.
The infinitesimal operator $v$ and its flow are given respectively
by changing the coefficient of $\frac{\partial}{\partial t}$ to
$t$ in (\ref{eq:19}) and the component $t+s$ of right hand side of
(\ref{eq:20}) to $t\,e^s$.\par
By fixing $(t,{\bf x},{\bf p})$ and acting the symmetry group on
it, then the flow of $v$, under the projection to the coordinates
in terms of $t$, $x_i$ and $p_j$ for a $1\leq j\leq n$, will be
similar to the space curve $(s,e^s,s^2+s)$ and has two branches,
the first of space exponential map and the second of a space
parabolic.
\paragraph{Example 3.6} For the case which we assume that $T=c\,t$ when $c$ is a
constant, and $K_i=1$ when $i=1,\cdots,n$, the contact
infinitesimal operator $v$ is as follows
\begin{eqnarray*}
v = c\,t\frac{\partial}{\partial t}+
\sum_i\,\Big\{\frac{\partial}{\partial
x_i}+\Big(\sum_{j,k}\,(a^i_{jk}+a^i_{kj})x_k\Big)
\frac{\partial}{\partial p_i}\Big\} \label{eq:21}
\end{eqnarray*}
The one--parameter of $v$ then is
\begin{eqnarray*}
\hspace{-2cm}(t,{\bf x},{\bf p}) \!\!\!\!& \mapsto&\!\!\!\!
\Big(t\,e^{cs}, {\bf x}+s, p_m + s^2\sum_l\, a^m_{ll}+ s
\sum_{l,k}\,(a^m_{lk}+a^m_{kl})x_k\Big)\label{eq:22}
\end{eqnarray*}
when in it, $m$ varies from 1 to $n$ and we assumed that $1\leq
k,l\leq n$. The action of the symmetry group on a fixed point,
will give a flow that its projection to three arbitrary axis
$(t,x_i,p_j)$ is similar to the 3-dimensional flow explained in
example 3.5. In the particular case of this example, the
characteristic equations tend to the following independent
differential invariants for $1\leq\alpha\leq n$ and
$n+1\leq\beta\leq 2n$:
\begin{eqnarray*}
\hspace{-0.3cm}&I_{\alpha}(t,{\bf x},{\bf p}) =
t(1-c\,k_{\alpha}), \hspace{0.6cm}I_{\beta}(t,{\bf x},{\bf p}) =
t\,\Big(\sum_{j,k}(a^{\beta-n}_{jk}+a^{\beta-n}_{kj})x_k\Big) -
c\,t\,p_{\beta-n}.&
\end{eqnarray*}
\paragraph{Example 3.7} By assuming $T=0$ and $K_i=x_i$ we find
that the infinitesimal operator as
\begin{eqnarray*}
v=\sum_i\,\Big\{x_i\frac{\partial}{\partial
x_i}+\Big(\sum_{j,k}\,(a^i_{jk}+a^i_{kj})x_j\,x_k\Big)
\frac{\partial}{\partial p_i}\Big\}
\end{eqnarray*}
and the one--parameter group of $v$ as
\begin{eqnarray*}
(t,{\bf x}, {\bf p})&\mapsto&\Big(t,e^s{\bf x}, p_i +
(e^{2s}-1)\sum_{j,k}\,(a^i_{jk}+a^i_{kj})x_j\,x_k\Big)
\end{eqnarray*}
where $i$ varies on values $1,\cdots,n$. Its flow also signifies a
parabolic if we reduce the jet coordinate $(t,{\bf x},{\bf p})$ to
$(t,x_l,p_k)$ for some specified $l$ and $k$. The concluded
invariants then will be as follows ($1\leq\alpha,\beta\leq n$)
\begin{eqnarray*}
&&\hspace{-0.5cm} I_1=t,\hspace{2cm} I_{\alpha}=\ln x_{\alpha}-\ln x_{\alpha-1},\\
&&\hspace{-0.5cm}
I_{\beta+n}=\sum_{j,k\neq\beta}\,(a^{\beta}_{jk}+a^{\beta}_{kj})x_j\,x_k\,x_{\beta}+
\sum_{k\neq\beta}\,(a^{\beta}_{{\beta}k}+a^{\beta}_{k{\beta}})x_k\,x_{\beta}+
a^{\beta}_{{\beta}{\beta}}\,x_{\beta}^2 - p_{\beta}.
\end{eqnarray*}
Finally, we present an example of the case in which $K_i$s are
functions dependent to $p_j$s resp.
\paragraph{Example 3.8} As we indicate, we suppose that $T=0$ and
$K_i=p_i$ for $i=1,\cdots,n$. So the infinitesimal generator is in
the form of
\begin{eqnarray*}
v = \sum_{i}\,\Big\{p_i\frac{\partial}{\partial x_i}+
\Big(\sum_{j,k}\,p_j(a^i_{jk}+a^i_{kj})x_k\Big)\frac{\partial}{\partial
p_i}\Big\}, \label{eq:19}
\end{eqnarray*}
that has the following trajectory in term of the parameter $s$
\begin{eqnarray*}
(t,{\bf x}, {\bf p})&\mapsto&(t,{\bf x}+{\bf p}s,{\bf p'})
\end{eqnarray*}
where ${\bf p'}$ consists of $n$ component, which its $i^{th}$
component is
\begin{eqnarray*}
p'_i= p_i +
\Big(\exp\Big((\sum_{k}(a^i_{ik}+a^i_{ki})x_k)s\Big)-1\Big)
\left(\frac{\sum_{j,k}\,p_j(a^i_{jk}+a^i_{kj})x_k}{\sum_{k}(a^i_{ik}+a^i_{ki})x_k}
\right).
\end{eqnarray*}
One can derive the following independent differential invariants
\begin{eqnarray*}
&&\hspace{-0.7cm} I_1=t,\hspace{2cm} I_{\alpha}=p_{\alpha}x_{\alpha-1}-p_{\alpha-1}x_{\alpha},\\
&&\hspace{-0.7cm}
I_{\beta+n}=\frac{p_{\beta}}{\sum_k\,(a^{\beta}_{\beta k} +
a^{\beta}_{k\beta}) x_k}-
\ln\Big(\sum_{j,k}\,p_j(a^{\beta}_{jk}+a^{\beta}_{kj})x_k\Big)
\frac{\sum_{j\neq\beta,k}\,p_j(a^{\beta}_{jk} +
a^{\beta}_{kj})x_k}{(\sum_k\,(a^{\beta}_{\beta
k}+a^{\beta}_{k\beta})x_k)^2} ,
\end{eqnarray*}
for $\alpha,\beta=1,\cdots,n$. Clearly, projection of the jet
coordinate to the coordinate $(x_m,p_l)$ for specified amounts $m$
and $l$, is the graph of the exponential map.\\

According to \cite{To}, if $G$ acts smoothly and transitively on a
manifold $M$, then $M$ is isomorphic to $G/H$ as a homogeneous
space, which obtained by quotienting by a closed Lie subgroup that
is, an isotropy subgroup provided by element of $G$. Thus, we can
use of the exponential map, $\exp:TM\rightarrow M$, defined in a
small neighborhood of identity element of $M$ as a Lie group. From
Baker--Campbell--Hausdorff formula we know that for vector fields
$X,Y$ of the Lie algebra of a Lie group, if we have $[X,Y]=0$,
then the flow of $X+Y$ is
$\xi_t=\phi_t\circ\psi_t=\psi_t\circ\phi_t$ where $\phi_t$ and
$\psi_t$ are flows of resp. $X$ and $Y$.
\paragraph{Remark 3.8} {\it With similar method as indicated in
above examples, by different selections of coefficients of
(\ref{eq:15}) as the linear combinations of 1-jet space variables,
their corresponding invariants are in the forms of (\ref{eq:18}),
and their proper projections of trajectories to 2 or 3-dimensional
space of their variables, providing their commutator be zero, are
similar to parabolic, exponential map or compositions of two
branches in the similar shape to them.}\\

This result is inferred, since the symmetric group prepared in the
above examples act transitively, if we add the condition that the
Lie bracket of every two infinitesimal operators be zero, then the
flow of every linear combination of any two of infinitesimal
generators is a composition of flows made by these two operators.


\begin{thebibliography}{}
\bibitem{CJ} C. Chicone and T. Jinghuang, {\em On General Proprties of Quaderatic Systems},
Amer. Math. Month. Vol. 89 {\bf 3}, 167--178 (1982).
\bibitem{Ha} Harrison, {\em The Differential Form Method for Finding Symmetries},
SIGMA {\bf 1}, 1--12 (2005).
\bibitem{Ne} T.A. Newton, {\em Two Dimensional quadratic homogeneous differential systems},
SIAM Rev., Vol. 20, {\bf 1} (1978).
\bibitem{Ib} N.H. Ibragimov, {\em Elementary Lie Group Analysis and Ordinary Differentail Equations},
John Wiley \& Sons, England (1999).
\bibitem{KLR} A. Kushner, V. Lychagin and V. Rubtsov, {\em Contact Geometry and
Non-linear Differential Equations}, Cambridge Univ. Press,
Cambridge (2007).
\bibitem{Ol} P.J. Olver, {\em Equivalence, Invariants, and
Symmetry}, Cambridge Univ. Press, Cambridge (1995).
\bibitem{To} Ph. Tondeur{\em Introduction to Lie groups and transformation groups},
Lecture Notes in Math. 7, Springer-Verlag, Berlin-Heidelberg-New
York, (1966).
\bibitem{Ts} A. Tsyvintsev, {\em On the existence of polynomial first integrals of
quadratic homogeneous systems of ordinary differential equations},
J. Phys. A: Math. Gen. {\bf 34} 2185–-2193 (2001).

\end{thebibliography}
\end{document}